\documentstyle[amssymb,amsfonts,12pt]{amsart}

\newtheorem{theorem}{Theorem}[section]

\def\C{{\mbox{\rm\kern.24em
\vrule width.03em height1.43ex depth-.052ex \kern-.26em C}}}
\def\QSet{\mbox{\rm\kern.24em
\vrule width.03em height1.48ex depth-.051ex \kern-.26em Q}}
\def\Z{{\mbox{\rm\kern.25em
\vrule width.03em height0.57ex depth0ex
\kern.033em
\vrule width.03em height1.52ex depth-0.96ex \kern-.338em Z}}}
\def\R{{\mbox{\rm I\kern-.22em R}}}

\def\sgn{{\rm sgn}}

\def\be#1{\begin{equation}\label{#1}}
\def\bas{\begin{align*}}
\def\eas{\end{align*}}
\def\bi{\begin{itemize}}
\def\ei{\end{itemize}}

\def \endprf{\hfill  {\vrule height6pt width6pt depth0pt}\medskip}
\def\emph#1{{\it #1}}

\title{A counterexample to a multilinear endpoint question of 
Christ and Kiselev}

\author{Camil Muscalu}
\address{Department of Mathematics, UCLA, Los Angeles CA 90095-1555}
\email{camil@@math.ucla.edu}

\author{Terence Tao}
\address{Department of Mathematics, UCLA, Los Angeles CA 90095-1555}
\email{tao@@math.ucla.edu}

\author{Christoph Thiele}
\address{Department of Mathematics, UCLA, Los Angeles CA 90095-1555}
\email{thiele@@math.ucla.edu}

\begin{document}

\begin{abstract}  
Christ and Kiselev \cite{ck0},\cite{ck1} have established that the generalized eigenfunctions of one-dimensional Dirac operators with $L^p$ potential $F$ 
are bounded for almost all energies for $p < 2$.  Roughly speaking, the proof involved writing these eigenfunctions as a multilinear series $\sum_n T_n(F, \ldots, F)$ and carefully bounding each term $T_n(F, \ldots, F)$.  It is conjectured that the results in \cite{ck1} also hold for $L^2$ potentials $F$.  However in this note we show that the bilinear term $T_2(F,F)$ and the trilinear term $T_3(F,F,F)$ are badly behaved on $L^2$, which seems to indicate that multilinear expansions are not the right tool for tackling this endpoint case.
\end{abstract}

\maketitle

\section{Introduction}

Let $F(x)$ be a real potential on $\R$.  For each energy $k^2 > 0$ we can consider the Dirac generalized eigenfunction equation 
$$ (\frac{d}{dx} + F) (- \frac{d}{dx} + F) \phi(x) = k^2 \phi(x)$$
on $\R$. This Dirac equation can be thought of as a Schr\"odinger equation with potential $V = F' + F^2$.  For each $k$ there are two linearly independent eigenfunctions $\phi = \phi_k$.  A natural question from spectral theory is to ask whether these eigenfunctions are bounded (i.e. are in $L^\infty_x$) for almost every real $k$.
In \cite{ck1} Christ and Kiselev\footnote{The results cited are phrased for Schr\"odinger operators but also extend to the slightly simpler case of Dirac operators, see \cite{ck3}.} showed among other things that this was true when $F \in L^p_x$ for any $1 \leq p < 2$.  It is well known (see e.g. \cite{pearson}) that the statement fails when $p > 2$, but the $p=2$ case remains open.  In \cite{deift} it is shown that for $L^2$ potentials one has absolutely continuous spectrum on $[0,\infty)$, but this is a slightly weaker statement.

We briefly outline the arguments in \cite{ck0},\cite{ck1}.  The method of variation of constants suggests the ansatz
\bas
\phi(x) &= a(x) e^{ikx} + b(x) e^{-ikx} \\
(-\frac{d}{dx} + F) \phi(x) &= -ik a(x) e^{ikx} + ik b(x) e^{-ikx}.
\end{align*}
Substituting this into the previous and simplifying, we reduce to the first-order system
\bas
a'(x) &= F(x) e^{-2ikx} b(x)\\
b'(x) &= F(x) e^{2ikx} a(x).
\end{align*}
For simplicity we may assume $F$ is supported on the positive half axis.
If we set initial conditions $a(-\infty) = 1$, $b(-\infty) = 0$ for instance, and then solve this system by iteration, we thus obtain the formal multilinear expansions
$$ a = 1 + \sum_{n \geq 2, \hbox{ even}} T_n(F, \ldots, F); \quad
b = \sum_{n \geq 1, \hbox{ odd}} T_n(F, \ldots, F)$$
where for each $n \geq 1$, $T_n$ is the $n$-linear operator 
$$ T_n(F_1, \ldots, F_n)(k,x) :=
\int_{x_1 < \ldots < x_n < x} e^{-2ik \sum_{j=1}^n (-1)^j x_j} F_1(x_1) \ldots F_n(x_n)\ dx_1 \ldots dx_n.$$
For integrable $F_j$ we can define the $n$-linear operators
$$ T_n(F_1, \ldots, F_n)(k,+\infty) :=
\int_{x_1 < \ldots < x_n} e^{-2ik \sum_{j=1}^n (-1)^j x_j} F_1(x_1) \ldots F_n(x_n)\ dx_1 \ldots dx_n.$$

The strategy of Christ and Kiselev was then to control each individual expression $T_n$ on $L^p$.  Specifally, they showed the estimate
\be{max}
\| \sup_x |T_n(F, \ldots, F)(k,x)| \|_{L^{p'/n,\infty}_k} \leq C_{p,n} \| F \|_{L^p_x}^n
\end{equation}
for all $n \geq 1$ and $1 \leq p < 2$, where $C_{p,n}$ was a constant which decayed rapidly in $n$ and $1/p + 1/p' := 1$.  In particular one has the non-maximal variant
\be{non-max}
\| T_n(F, \ldots, F)(k,+\infty) \|_{L^{p'/n,\infty}_k} \leq C_{p,n} \| F \|_{L^p_x}^n.
\end{equation}
The boundedness of eigenfunctions for almost every $k$ then follows by summing these bounds carefully.

It is tempting to try this approach for the endpoint $p=2$.  For $n=1$ we see that $T_1(F)(k,+\infty)$ is essentially the Fourier co-efficient $\hat F(k)$, while $\sup_x |T_1(F)(k,x)|$ is essentially the Carleson maximal operator $CF(k)$.  The estimates \eqref{non-max},  \eqref{max} for $p=2$ then follow from Plancherel's theorem and the Carleson-Hunt theorem \cite{carleson}, \cite{hunt} respectively.

For $n=2$ the expression $T_2(F,F)(k,+\infty)$ is essentially $H_-(|\hat F|^2)(k)$, where $H_-$ is the Riesz projection
$$ \widehat{H_- F} := \chi_{(-\infty,0]} \widehat{F},$$
and so \eqref{non-max} follows for $p=2$ by H\"older's inequality and the weak-type $(1,1)$ of the Riesz projections.  We also remark that if the phase function $x_1 - x_2$ in the definition of $T_2$ were replaced by $\alpha_1 x_1 + \alpha_2 x_2$ for generic numbers $\alpha_1$, $\alpha_2$ then the operator is essentially a bilinear Hilbert transform and one still has boundedness 
from the results in \cite{lt0}, \cite{lt1}, \cite{thiele}.

It may thus appear encouraging to try to estimate the higher order multilinear operators for $L^2$ potentials $F$.  However, in this note we show

\begin{theorem}\label{main}  When $p=2$ and $n=2$, the estimate \eqref{max} fails.  When $p=2$ and $n=3$, the estimate \eqref{non-max} fails.
\end{theorem}

Because of this, we believe that it is not possible to prove the almost
everywhere boundedness of eigenfunctions for Dirac or Schr\"odinger operators with $L^2$ potential purely by multilinear expansions; we discuss this further in the remarks section.

The counterexample has a logarithmic divergence, and essentially relies on the fact that while convolution with the Hilbert kernel $p.v. \frac{1}{x}$ is bounded, convolution with $\frac{\sgn(x)}{x}$ or $\frac{\chi_{(-\infty,0]}(x)}{x}$ is not.  It may be viewed as an assertion that $L^2$ potentials create significant long-range interaction effects which are not present for more rapidly decaying potentials.

Interestingly, our counterexamples rely strongly on a certain degeneracy in the phase function $\sum_j (-1)^j x_j$ on the boundary of the simplex $x_1 < \ldots < x_n$.  If one replaced this phase by $\sum_j x_j$, then we have shown in \cite{mtt0}, \cite{mtt1} that the bound \eqref{non-max} in fact holds when $p=2$ and $n=3$.  Indeed this statement is true for generic phases of the form $\sum_j \alpha_j x_j$.  A similar statement holds for \eqref{max} when $p=2$ and $n=2$ and will appear elsewhere.  

The first author was supported by NSF grant DMS 0100796.
The second author is a Clay Prize Fellow and is supported by a grant from
the Packard Foundation. The third author was supported
by a Sloan Fellowship and by NSF grants DMS 9985572 and DMS 9970469.
The authors are grateful to M. Christ for pointing out the importance of the degeneracy in the phase function $\sum_j (-1)^j x_j$ and for suggesting numerous valuable improvements to the manuscript.

\section{Proof of Theorem \ref{main}}

The letter $C$ may denote different large constants in the sequel.  To be consistent with the previous notation we shall define the Fourier transform as
$$ \hat F(k) := \int e^{-2ikx} F(x)\ dx.$$

We let $N \gg 1$ be a large integer parameter, which we shall take to be a square number, and test \eqref{max}, \eqref{non-max} with the real-valued potential
$$ F(x) := \sum_{j=N}^{2N} F_j(x)$$
where the $F_j$ are given by
$$ F_j(x) := N^{-1} \cos(2 \frac{Aj}{N} x) \phi(\frac{x}{N} - j),$$
$\phi$ is a smooth real valued function supported in $[-\frac{1}{4},\frac{1}{4}]$ with total mass $\int \phi = 1$ such that $\widehat{\phi}$
stays away from $0$ in $[-1,1]$, and
$A$ is a sufficiently large absolute constant whose purpose is to ensure
that 
$$4 \sum_{j\in \Z\setminus\{0\}} \left| \widehat {\phi}(\xi-Aj) \right|
\le |\widehat{\phi(\xi)}|$$
for $\xi \in [-1,1]$.
Informally, $F$ is a ``chirp'' which is localized in phase space to the region 
$$\{ (k,x): k = \pm \frac{Aj}{N} + O(\frac{1}{N}); x = Nj + O(N), N \leq j \leq 2N \}.$$

We may compute the Fourier transform of the $F_j$ using the rapid decay of $\hat \phi$ as
\be{fj-ft}
\hat F_j(k) = \frac{1}{2} e^{-2i(Nk-Aj)j} \hat \phi(Nk - Aj) + O(N^{-200})
\end{equation}
in the region $\frac{A}{2} < k < 3A$.  We remark that the error term $O(N^{-200})$ has a gradient which is also $O(N^{-200})$.

Clearly we have $\|F_j\|_2=O(N^{-1/2})$, and hence that
$$ \| F \|_2 = O(1).$$
We now compute 
\be{2-def}
T_2(F, F)(k,x) = \int_{x_1 < x_2 < x} e^{ 2ik(x_1-x_2)} F(x_1) F(x_2)\ dx_1 dx_2
\end{equation} 
in the region
\be{j-region}
\left|Nk - {Aj_0}\right|\le  1; \quad x = N(j_0-\sqrt{N} + \frac{1}{2}) 
\end{equation}
for some integer $\frac{3N}{2} < j_0 < 2N$.  In this region we show that
\be{2-bound}
|T_2(F,F)(k,x)| \geq C^{-1} \log N,
\end{equation}
which will imply that
$$ \| \sup_x |T_2(F,F)|(k,x) \|_{L^{2,\infty}_k} \geq C^{-1} \log N$$
and thus contradict \eqref{max} for $n=2$ and $p=2$ by letting $N$ go to infinity.

We now prove \eqref{2-bound}.  Fix $k$, $j_0$, $x$.  Observe from \eqref{2-def} that $T_2(F_j,F_{j'})(k,x)$ vanishes unless $j \leq j' \leq j_0 - \sqrt{N}$.  Thus we may expand
\begin{align}
T_2(F,F)(k,x) &=  \sum_{N \leq j \leq j_0 - \sqrt{N}} T_2(F_j,F_j)(k,x) \label{jj}\\
&+ \sum_{N \leq j < j' \leq j_0 - \sqrt{N}} T(F_j,F_{j'})(k,x) \label{jjp}.
\end{align}
We first dispose of the error term \eqref{jjp}.  In the region $j < j' \leq j_0 - \sqrt{N}$, the conditions $x_1 < x_2 < x$ in \eqref{2-def} become superfluous, so we may factor
$$ T_2(F_j,F_{j'})(k,x) = \overline{\hat{F_j}(k)} {\hat{F_{j'}}(k)}.$$
However, since $\hat \phi$ is rapidly decreasing and $|j-j_0|, |j'-j_0| \geq \sqrt{N}$, we see from \eqref{fj-ft} that
$$ |\hat{F_j}(k)|, |\hat{F_{j'}}(k)| \leq C N^{-100}.$$
Summing this, we see that the total contribution of \eqref{jjp} is $O( N^{-198})$.

Now we consider the contribution of \eqref{jj}.  We use the identity
\be{t2}
T_2(F_j,F_j)(k,x) = T_2(F_j,F_j)(k,+\infty)=H_-( |\hat F_j|^2 )(k)
\end{equation} 
combined with \eqref{fj-ft}.  The operator $H_-$ is a non-trivial linear combination of the identity and the Hilbert transform, while $|\hat F_j|^2$ is essentially a non-negative bump function rapidly decreasing away from the interval $[jA/N-O(1/N), jA/N+O(1/N)]$.  Because of this we see that for $j\neq j_0$ we
have
\be{t2-j}
H_-( |\hat F_j|^2 )(k) =  \frac{c}{j-j_0} + O(|j-j_0|^{-2})
\end{equation}
where $c$ is a non-zero absolute constant.  Summing this over all $j \leq j_0 - \sqrt{N}$  and observing that $j - j_0$ has a consistent sign we see that the contribution of \eqref{jj} has magnitude at least
$C^{-1} \log N$, and \eqref{2-bound} follows.

We now compute $T_3(F, F, F)(k,+\infty)$ in the region
\be{j-region-2}
\left|Nk - {Aj_0}\right|\le 1; \quad 1.4 N < j_0 < 1.6 N.
\end{equation}
We will show that
\be{3-bound}
|T_3(F,F,F)(k,+\infty)| \geq C^{-1} \log N
\end{equation}
in this region, which will disprove \eqref{non-max} for $n=3$ and $p=2$ similarly to before.

It remains to prove \eqref{3-bound}.  Fix $j_0$.  Observe that $T_3(F_j,F_{j'},F_{j''})(k,+\infty)$ vanishes unless $j \leq j' \leq j''$.  Thus we can split
\begin{align}
T_3(F,F,F)(k,+\infty) &=  \sum_{N \leq j \leq 2N} T_3(F_j,F_j,F_j)(k,+\infty) \label{kkk}\\
&+ \sum_{N \leq j < j' \leq 2N} T_3(F_j,F_j,F_{j'})(k, +\infty) \label{kkp}\\
&+ \sum_{N \leq j' < j \leq 2N} T_3(F_{j'}, F_j, F_j) (k, +\infty) \label{kpp}\\
&+ \sum_{N \leq j < j' < j'' \leq 2N} T_3(F_j, F_{j'}, F_{j''})(k, +\infty). \label{klm}
\end{align}

We first consider \eqref{kkk}.  We expand
$$ T_3(F_j,F_j,F_j)(k, +\infty) = \int_{x_1 < x_2 < x_3} e^{2ik(x_1-x_2+x_3)} F_j(x_1) F_j(x_2) F_j(x_3)\ dx_1 dx_2 dx_3.$$
This is a linear combination of eight terms of the form
$$ N^{-3} \int_{x_1 < x_2 < x_3} e^{2ik(x_1-x_2+x_3)} 
e^{2i\frac{Aj}{N}(\pm x_1 \pm x_2 \pm x_3)}
\phi(\frac{x_1}{N} - j) \phi(\frac{x_2}{N} - j) \phi(\frac{x_3}{N} - j) 
\ dx_1 dx_2 dx_3;$$
making the substitutions $y_s := \frac{x_s}{N} - j$ for $s=1,2,3$, this becomes
$$ e^{i\theta} \int_{y_1 < y_2 < y_3} e^{2ikN(y_1-y_2+y_3)} 
e^{2i Aj (\pm y_1 \pm y_2 \pm y_3)}
\phi(y_1) \phi(y_2) \phi(y_3) 
\ dy_1 dy_2 dy_3$$
for some phase $e^{i\theta}$ depending on all the above variables.

We shall only consider the choice of signs $(- y_1 + y_2 - y_3)$; the reader may easily verify that the other choices of signs are much smaller thanks to stationary phase.  In this case we can write the above as
$$ e^{i\theta} \int_{y_1 < y_2 < y_3} e^{2i(kN - Aj)(y_1-y_2+y_3)} 
\phi(y_1) \phi(y_2) \phi(y_3) 
\ dy_1 dy_2 dy_3.$$
If $kN - Aj = O(1)$ we estimate this crudely by $O(1)$.  Otherwise we can perform the $y_1$ integral using stationary phase to obtain
$$ e^{i\theta} \frac{1}{2i(kN-Aj)} \int_{y_2 < y_3} e^{2i(kN - Aj) y_3} 
\phi(y_2) \phi(y_2) \phi(y_3) 
\ dy_2 dy_3 + O(|kN-Aj|^{-2}).$$
Performing another stationary phase we see that this quantity is
$O(|kN-Aj|^{-2})$.
Summing over all $j$ we see that \eqref{kkk} is $O(1)$.

Let us now consider \eqref{klm}.  When $j < j' < j''$, the constraints $x_1 < x_2 < x_3$ in the definition of $T_3$ are redundant, and we can factorize
$$ T_3(F_j, F_{j'}, F_{j''})(k, +\infty)
= \overline{\hat{F_j}(k)} {\hat F_{j'}(k)} \overline{\hat F_{j''}(k)}.$$
Applying \eqref{fj-ft} and using the rapid decay of $\hat \phi$ we see that
$$ | T_3(F_j, F_{j'}, F_{j''})(k, +\infty) |
\leq C (1 + |j - j_0| + |j' - j_0| + |j'' - j_0|)^{-10} + CN^{-100}.$$
Summing over all $j,j',j''$ we see that \eqref{klm} is $O(1)$.

It remains to control $\eqref{kpp} + \eqref{kkp}$.  First we consider \eqref{kkp}. For this term the condition 
$x_2 < x_3$ is redundant, so we can factorize
$$ T_3(F_j,F_j,F_{j'})(k, +\infty) = T_2(F_j,F_j)(k, +\infty) 
\overline{\hat F_{j'}(k)}.$$

Now consider \eqref{kpp}.  For this term the condition $x_1 < x_2$ is redundant, so we can factorize
$$ T_3(F_{j'},F_j,F_j)(k, +\infty) = \overline{\hat F_{j'} (k)}
\int_{x_2 < x_3} e^{2ik(x_3-x_2)} F_j(x_2) F_j(x_3)\ dx_3 dx_2.
$$
Writing $x_1$ instead of $x_3$ we thus have
$$ T_3(F_{j'},F_j,F_j)(k, +\infty) = \overline{\hat F_{j'} (k)}
(|\hat F_j(k)|^2 - T_2(F_j,F_j)(k, +\infty)).$$
Combining this with the previous, we thus see that
$$ \eqref{kpp} + \eqref{kkp} =
\sum_{N \leq j, j' \leq 2N} \sgn(j'-j) T_2(F_j,F_j)(k,+\infty) 
\overline{\hat F_{j'}(k)}
+ \sum_{N \leq j' < j \leq 2N} 
\overline{\hat F_{j'} (k)} |\hat F_j(k)|^2.$$ 
Using \eqref{fj-ft} as in \eqref{klm} we see the second term is $O(1)$, so to prove \eqref{3-bound} it will suffice to show that
\begin{equation}\label{oneovermodx}
|\sum_{N \leq j, j' \leq 2N, j\neq j'} \sgn(j'-j) T_2(F_j,F_j)(k,+\infty) 
\overline{\hat F_{j'}(k)}|
\geq C^{-1} \log N.
\end{equation}
We first consider the terms with $j'=j_0$. We claim these terms are
the dominant contribution. 
From \eqref{t2}, \eqref{t2-j} we conclude
$$
\sum_{N \leq j \leq 2N, j\neq j_0} \sgn(j_0-j) T_2(F_j,F_j)(k,+\infty) 
\overline{\hat F_{j_0}(k)}$$
\begin{equation}\label{dominant}=
\sum_{N \leq j \leq 2N, j\neq j_0} c \frac{\sgn(j_0-j)}{j_0-j}
\overline{\hat F_{j_0}(k)}\ +\ O(1)\ \ .
\end{equation}
Here $c$ is the same non-zero constant as in \eqref{t2-j}, and
$\hat F_{j_0}(k) $ is bounded away from $0$ by choice of $\phi$.
Thus the first term is greater than $C^{-1}\log N$, so it suffices indeed
to show that this term is the dominant contribution to (\ref{oneovermodx}).

We consider the terms with $j=j_0$. Using that
$|T_2(F_j,F_j)(k,\infty)|\le C$ we obtain 
$$
\sum_{N \leq j' \leq 2N, j_0\neq j'} |T_2(F_{j_0},F_{j_0})(k,+\infty) 
\overline{\hat F_{j'}(k)}|
\le C$$
This term is therefore negligible.

Finally, we have to consider the terms with $j,j'\neq j_0$.
We have by the choice of $A$,
$$
\sum_{N \leq j,j' \leq 2N, j,j'\neq j_0} |T_2(F_j,F_j)(k,+\infty)| 
|\overline{\hat F_{j'}(k)}|$$
$$
\le \frac 12 \sum_{N \leq j \leq 2N, j\neq j_0} \frac c{|j-j_0|} 
|\overline{\hat F_{j_0}(k)}|+ C$$
This term is dominated by (\ref{dominant}).
This completes the proof of \eqref{3-bound}.

\section{Remarks}\label{remarks-sec}

\begin{itemize}
\item The counterexample can easily be extended to larger $n$ (e.g. by appending some bump functions to the left or right of $F$).

\item The counterexample above involved a potential $F$ which was bounded in $L^2$, but for which $\sup_x |T_2(F,F)(k,x)|$ and $|T_3(F,F,F)(k,+\infty)|$ were large (about $\log N$) on a large subset of $[A,2A]$.  By letting $N$ vary and taking suitable linear combinations of such variants of the above counterexample, one can in fact generate a potential $F$ bounded in $L^2$ for which $\sup_x |T_2(F,F)(k,x)|$ is infinite and $|T_3(F,F,F)(k,x)|$ accumulates at $\infty$ for $x\to \infty$ for all $k$ in a set of positive measure (one can even achieve
blow-up almost everywhere). Thus it is not possible to estimate these multilinear expansions in any reasonable norm if one only assumes the potential to be in $L^2$.  
Similarly if $F$ had a derivative
in $L^2$; it is the decay of $F$ which is relevant here, not the regularity.

\item The unboundedness of $T_3$ on $L^2$ can be interpreted as stating that the (non-linear) scattering map $F \mapsto b_k(+\infty)$ from potentials to reflection coefficients is not $C^3$ on the domain of $L^2$ potentials.  Similarly the map $F \mapsto a_k(+\infty)$ from potentials to transmission coefficients is not $C^4$ on the domain of $L^2$ potentials.  In particular these scattering maps are not analytic.

\item Despite the bad behavior of the individual terms $T_k(F, \ldots, F)$, the transmission and reflection coefficients $a_k(x)$, $b_k(x)$ are still bounded for the counterexample given above.  This phenomenon is similar to the observation that the function $e^{ix} = 1 + ix - x^2/2 - \ldots$ is bounded for arbitrarily large \emph{real} $x$, even if the individual terms $(ix)^n/n!$ are not.

We now sketch the proof of boundedness of $a_k$, $b_k$.  Suppose that $k = Aj_0/N + O(1/N)$ for some $N \leq j_0 \leq 2N$; we now fix $j_0$ and $k$.  We can write
$$ 
\left( 
\begin{array}{l}
a_k(x) \\
b_k(x)
\end{array}
\right)
=
G(x) 
\left( 
\begin{array}{l}
1 \\
0
\end{array}
\right)
$$
where $G$ is the $2 \times 2$ matrix solving the ODE
$$ G'(x) = 
\left( 
\begin{array}{ll}
0 & F(x) e^{-2ikx} \\
F(x) e^{2ikx} & 0
\end{array}
\right) G(x); \quad G(-\infty) =
\left( 
\begin{array}{ll}
1 & 0 \\
0 & 1
\end{array}
\right).
$$
We define the matrices $G_j$ similarly by
$$ G_j'(x) = 
\left( 
\begin{array}{ll}
0 & F_j(x) e^{-2ikx} \\
F_j(x) e^{2ikx} & 0
\end{array}
\right) G_j(x); \quad G_j(-\infty) =
\left( 
\begin{array}{ll}
1 & 0 \\
0 & 1
\end{array}
\right).
$$
We observe the identity
$$ G(x) = G_{j_1}(+\infty) G_{j_1-1}(+\infty) \ldots G_N(+\infty)$$
whenever $N \leq j_1 \leq 2N$ and $x = N (j_1 + \frac{1}{2})$; this can be proven by an easy induction on $j_1$ and the observation that the above ODE are invariant under right-multiplication.

One can compute the $G_j(+\infty)$ using multilinear expansions (or using Gronwall's inequality), eventually obtaining
$$ G_j(+\infty) = 
\left( 
\begin{array}{ll}
1 + \frac{iC}{j-j_0} & 0 \\
0 & 1 - \frac{iC}{j-j_0}
\end{array}
\right) + O(|j-j_0|^{-2})
$$
for all $j \neq j_0$, where $C$ is a non-zero real constant.  Because of the crucial factor of $i$ in the diagonal entries we see that the operator norm $\| G_j(+\infty)\|$ of $G_j$ is
$$ \| G_j(+\infty)\| = 1 + O(|j-j_0|^{-2}).$$
This allows one to multiply the $G_j(+\infty)$ together and obtain boundedness of $G(x)$ and hence $a_k(x)$, $b_k(x)$.

In analogy with the observation concerning $e^{ix}$, one may need to use the fact that $F$ is real in order to obtain boundedness of eigenfunctions in the $L^2$ case.  When $F$ is real there are additional estimates available, such as the scattering identity
$$ \int \log |a_k(+\infty)|\ dk = C \int |F(x)|^2\ dx$$
for some absolute constant $C$; see for instance \cite{deift}.

We do not yet know how to obtain boundedness of eigenfunctions for $L^2$ potentials $F$.  However we have been able to achieve this for a model problem in which the Fourier phases $e^{2ikx}$ are replaced by a dyadic Walsh variant $e(k, x)$.  See \cite{mtt:walsh}.

\item One can modify the counterexample to provide similar counterexamples for Schr\"odinger operators $-\frac{d^2}{dx^2} + V$ with $V\in L^2$, either by using the Miura transform $V = F' + F^2$ mentioned in the introduction, or by inserting the standard WKB phase modification to the operators $T_k$ as in \cite{ck0}.  We omit the details.

\item The multilinear expansion of $a$ leads to an expansion
of $|a|^2$, whose quadratic term is equal to
$$2 Re (T_2(F,F))= 2 Re (H_-(|\widehat{F}|^2)) = |\widehat{F}|^2$$
This term is in $L^1$, which is better than the term $T_2(F,F)$,
which is in general only in the Lorentz space $L^{1,\infty}$. The higher
order terms of the expansion of $|a|^2$ are however unbounded
again. Using the identity $|a|^2=1+|b|^2$ we see that the
fourth order term of $|a|^2$ is equal to
$$2 Re (\overline{T_1(F)} T_3(F,F,F))$$
We now define the modified potential
$$G(x)=F(x)+G_0(x)$$
where $F$ is as in the proof of Theorem \ref{main}
and $G_0(x)= \phi(x-N^3)$.
Expanding the fourth order term by multilinearity, one observes
that all terms can be estimated from above nicely with the exception of
$$2 Re (\overline{T_1(G)} T_3(F,F,F))$$
Since $T_1(G)=\widehat{G}$ has more rapidly changing phase
than $T_3(F,F,F)$, the real part and the modulus 
$\overline{T_1(G)} T_3(F,F,F)$ are of comparable size
on a large set, and so this term is of the order $\log(N)$ on 
a large set just like  $T_3(F,F,F)$ itself.

\end{itemize}


\begin{thebibliography}{99}

\bibitem{carleson} Carleson, L,
\emph{On convergence and growth of partial sums of Fourier series}
Acta Math. \textbf{116} [1966], pp. 135-157.
 
\bibitem{ck0} Christ, M., Kiselev, A.,
\emph{Maximal functions associated to filtrations}, J. Funct. Anal.
\textbf{179} [2001], no 2., pp. 409-425

\bibitem{ck1} Christ, M., Kiselev, A.,
\emph{WKB asymptotic behaviour of almost all 
generalized eigenfunctions of one-dimensional
Schr\"{o}dinger operators with slowly decaying potentials}, 
J. Funct. Anal.
\textbf{179} [2001], no 2., pp. 426-447

\bibitem{ck3} Christ, M., Kiselev, A.,
\emph{Scattering and wave operators for one-dimensional Schr\"odinger operators with nonsmooth slowly decaying potentials}, preprint.

\bibitem{deift} Deift, P. and Killip, R.,
\emph{On the absolutely continuous spectrum of one-dimensional Schr\"odinger operators with square summable potentials}, Comm. Math. Phys. \textbf{203} [1999], pp. 61--72.

\bibitem{hunt} Hunt, R. \emph{On the convergence of Fourier series}
\emph{1968 Orthogonal Expansions and their Continuous Analogues. Proc
Conf. Edwardsville, Ill. 1967} pp. 235-255, Southern Illinois Univ. Press,
Carbondale, Ill.

\bibitem{lt0} Lacey, M., Thiele, C., \emph{$L^p$ Bounds for the bilinear Hilbert 
transform, 
$p>2$.}  { Ann. Math.} (2) \textbf{146} [1997] no. 2, pp. 693---724.

\bibitem{lt1}  Lacey, M.,  Thiele, C.,  \emph{ On Calder\`on's Conjecture.}  
 { Ann.  Math.} (2) \textbf{149} [1999], pp. 475--496

\bibitem{mtt0} Muscalu, C., Tao, T., Thiele, C.,
\emph{$L^p$ estimates for the biest I. The Walsh case.} preprint

\bibitem{mtt1} Muscalu, C., Tao, T., Thiele, C.,
\emph{$L^p$ estimates for the biest II. The Fourier case.} preprint

\bibitem{mtt:walsh} Muscalu, C., Tao, T., Thiele, C.,
\emph{A Walsh model for the bounded eigenfunction problem for $L^2$ potentials}, preprint

\bibitem{pearson} Pearson, D., \emph{Singular continuous measures in scattering theory}, Comm. Math. Phys. \textbf{60} [1978], pp. 13--36.

\bibitem{thiele} Thiele, C. {\it On the Bilinear Hilbert transform}. 
Universit\"at Kiel, Habilitationsschrift [1998]


\end{thebibliography}
\end{document}